\documentclass[a4paper,12pt]{article}
\usepackage{mathrsfs}
\usepackage{}
\usepackage[top=2.0cm,bottom=2.0cm,left=2.0cm,right=2.0cm]{geometry}
\usepackage{amssymb}
\usepackage{graphicx}
\usepackage{amsmath,amsthm,amssymb,lineno}
\usepackage{latexsym}
\usepackage{epstopdf}
\usepackage{setspace}
\usepackage{amsmath,color}
\usepackage{graphicx,booktabs,multirow}
\usepackage{latexsym, tabularx,shapepar}
\usepackage[all,2cell,dvips]{xy} \UseAllTwocells \SilentMatrices
\usepackage{appendix}
\usepackage{longtable}
\usepackage{cite}
\usepackage{CJK}
\usepackage{indentfirst}
\usepackage{array}
\usepackage{amsmath}
\usepackage[colorlinks,
           linkcolor=blue,
           anchorcolor=blue,
           citecolor=blue
           ]{hyperref}
\usepackage{dsfont}

\usepackage{etoolbox}
\let\bbordermatrix\bordermatrix
\patchcmd{\bbordermatrix}{8.75}{4.75}{}{}
\patchcmd{\bbordermatrix}{\left(}{\left[}{}{}
\patchcmd{\bbordermatrix}{\right)}{\right]}{}{}

\allowdisplaybreaks[4]
\graphicspath{{figures/}}
\usepackage{caption}
\makeatletter
\newcommand{\rmnum}[1]{\romannumeral #1}
\newcommand{\Rmnum}[1]{\expandafter\@slowromancap\romannumeral #1@}
\makeatother

\newtheorem{theorem}{Theorem}[section]
\newtheorem{lemma}[theorem]{Lemma}
\newtheorem{remark}[theorem]{Remark}

\newtheorem{proposition}[theorem]{Proposition}

\begin{document}

\title{ Extremal results for $\mathcal{K}^-_{r + 1}$-free signed graphs}

\author{Zhuang Xiong\\
	{\small College of Mathematics and Statistics, Hunan Normal University, zhuangxiong@hunnu.edu.cn}\\
	\and
		Yaoping Hou\\
		{\small College of Mathematics and Statistics, Hunan Normal University, yphou@hunnu.edu.cn}}

\maketitle
\begin{abstract}
	This paper gives tight upper bounds on the number of edges and the index for $\mathcal{K}^-_{r + 1}$-free unbalanced signed graphs, where $\mathcal{K}^-_{r + 1}$ is the set of $r+1$-vertices unbalanced signed complete graphs. \\
	\indent We first prove that if $\Gamma$ is an $n$-vertices $\mathcal{K}^-_{r + 1}$-free unbalanced signed graph, then the number of edges of $\Gamma$ is $$e(\Gamma) \leq \frac{n(n-1)}{2} - (n - r ).$$
   	\indent Let $\Gamma_{1,r-2}$ be a signed graph obtained by adding one negative edge and  $r - 2$ positive edges between a vertex and an all-positive signed complete graph $K_{n - 1}$.
   	Secondly, we show that if $\Gamma$ is an $n$-vertices $\mathcal{K}^-_{r + 1}$-free unbalanced signed graph, then the index of $\Gamma$ is 
   	$$\lambda_{1}(\Gamma) \leq \lambda_{1}(\Gamma_{1,r-2}), $$ with equality holding if and only if $\Gamma$ is switching equivalent to $\Gamma_{1,r-2}$. \\
   	\indent It is shown that these results are significant in extremal graph theory. 
   	Because they can be regarded as extensions of Tur{\'a}n's Theorem [Math. Fiz. Lapok 48 (1941) 436--452] and spectral Tur{\'a}n problem [Linear Algebra Appl. 427 (2007) 183--189] on signed graphs, respectively. 
   	Furthermore, the second result partly resolves a recent open problem raised by Wang [arXiv preprint arXiv:2309.15434 (2023)].
\end{abstract}
\baselineskip=0.25in

\section{ Introduction}
\begin{figure}[htbp]\centering\hspace{0cm}
\scalebox{0.5}{\includegraphics[width=25cm]{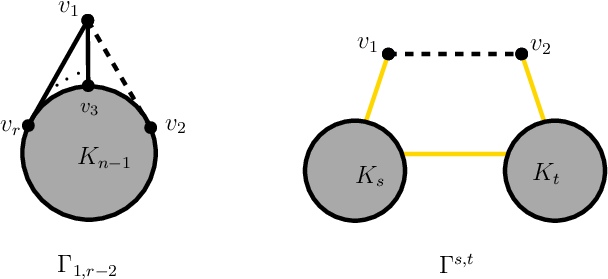}}\\
\caption{ The signed graphs mentioned in the introduction. } 
\label{fig:ExtreGs}
\end{figure}

  A signed graph $\Gamma$ of order $n$ is a pair $(G,\sigma)$, where $G = (V(G),E(G))$ is an $n$-vertices (unsigned) graph with vertex set $V(G)$ and edge set $E(G)$, called the underlying graph, and $\sigma: E(G) \rightarrow \{-1, +1\}$ is the sign function. 
  An edge $e = uv$ is positive (resp., negative) if $\sigma(e) =  +1$ (resp., $-1$) and denoted by $u \stackrel{+}{\sim} v$ (resp., $u \stackrel{-}{\sim} v$). 
  A signed graph $\Gamma$ is called homogeneous if all its edges have the same sign, and heterogeneous otherwise. If all edges are positive (resp., negative), then $\Gamma$ is called all-positive (resp., all-negative) and denoted by $(G, +)$ (resp., $(G, -)$). \\
  \indent A cycle in $\Gamma$ is called positive (resp., negative) if it contains even (resp., odd) negative edges. A signed graph is balanced if all its cycles are positive, otherwise it is unbalanced.
  The adjacency matrix of a signed graph $\Gamma$ is defined as $A(\Gamma) = (a_{ij}^{\sigma})$ , where $a_{ij}^{\sigma} = \sigma(v_iv_j)$ if $v_i \sim v_j$ and $a_{ij}^{\sigma} =0$ otherwise.  
  The eigenvalues of $\Gamma$ are identified with those of $A(\Gamma)$, and we denote them by $\lambda_1(\Gamma) \geq \lambda_2(\Gamma) \geq \dots \geq \lambda_n(\Gamma)$. The largest eigenvalue of $\Gamma$ is called the index of $\Gamma$ and the spectral radius of $\Gamma$ is defined as $\rho(\Gamma) = max\{\lambda_1(\Gamma), -\lambda_n(\Gamma) \}$.
  For $U \subset V(\Gamma)$, if we reverse the signs of all edges between $V(\Gamma)$ and $V(\Gamma) \setminus U$ to obtain a signed graph $\Gamma_U$, then $\Gamma$ is called switching equivalent to $\Gamma_U$, written as $\Gamma \sim \Gamma_U$. 
  Switching operation plays a key role in the study of spectral theory of signed graphs, because two switching equivalent signed graphs share the same spectrum. 
  In fact, any switching arising from $U$ can be realized by a diagonal matrix $S_U = diag(s_1, s_2, \dots, s_n)$ having $s_i = 1$ for each $i \in U$, and
  $s_i = -1$ otherwise. 
  Hence, $A(\Gamma) = S_U^{-1}A(\Gamma_U) S_U$; in this case we say that the matrices are signature
  similar.
  Throughout this paper we consider only simple signed graphs $\Gamma$ with $n$ vertices and $e(\Gamma)$ edges.
  For more notations and notions of signed graphs, we refer to 	\cite{zaslavskyi1982signed}.\\
  \indent Before presenting new theorems, we provide an introductory discussion. 
  Given a graph $F$, a graph $G$ is called $F$-free, if it contains no $F$ as a subgraph. 
  As the beginning of the extremal graph theory, in 1940, Tur{\'a}n raised and solved the extremal problem for $K_{r + 1}$ \cite{turan1941an,turan1954the}. 
  	Tur{\'a}n's graph, denoted by $T_r(n)$, is the complete $r$-partite graph on $n$ vertices which is the
  	result of partitioning $n$ vertices into $r$ almost equally sized partitions $(\lfloor n/r \rfloor , \lceil n/r \rceil) $ and taking all edges connecting two different partition classes (note that if $n \leq r$ then $T_r(n) = K_n$). Denote the number of
  	edges in Tur{\'a}n's graph by $t_r(n) = |E(T_r(n))|$.
  	Using these notations, Tur{\'a}n's Theorem can be stated as follows. 
  	\begin{theorem}\label{tho:turan}
  		(Tur{\'a}n's Theorem). If $G$ is $K_{r+1}$-free then $e(G) \leq t_r(n)$. Furthermore, equality holds if and only if $G = T_r(n)$.
  	\end{theorem}
	Subsequently, determining the maximum number of edges (referred to as the Tur{\'a}n number) in $n$-vertices $F$-free graphs became the central problem of classical extremal graph theory, known as the Tur{\'a}n problem.
	The Tur{\'a}n problem has been extensively studied, and readers can refer to \cite{bollobas1978extremal} for a comprehensive overview.\\
	\indent In the past two decades, the spectral version of Tur{\'a}n problem is paid much attention by many researchers.
	Below we only mention the result raised by Nikiforov in 2007\cite{nikiforovi2007bounds}, which will be used for study later in this article. 
	The readers can refer to an outstanding paper \cite{li2022survey} to understand the current research status of such problems. 
	\begin{theorem}\label{tho:nik'sgraph}\cite[Theorem 1]{nikiforovi2007bounds}
		If $G$ is $K_{r + 1}$-free then $\rho(G) \leq \rho(T_r(n))$. Furthermore, equality holds if and only if $G = T_r(n)$.
	\end{theorem}
  	Note that the edge extremal graph of Theorem \ref{tho:turan} is the same as the spectral extremal graph of Theorem \ref{tho:nik'sgraph}.
  	Actually, the classical Tur{\'a}n problem and the spectral Tur{\'a}n problem are closely related. See \cite{nikiforovi2011some} for more discussions on this topic.\\  
 	\indent Different from aforementioned studies, we focus on signed graphs, and ask what are the maximum number of edges and the maximum spectral radius of an $\mathcal{F}$-free signed graph of order $n$, where $\mathcal{F}$ is the set of signed graphs. 
	These problems were initially proposed by Wang, Hou, and Li in their recent research \cite{wang2022extremed}, in which they referred to these problem as the Tur{\'a}n-like problem in the context of signed graphs.
	Here we refer to these two problems as the signed Tur{\'a}n problem and signed spectral Tur{\'a}n problem. 
	Denote by $\mathcal{K}_k^-$ and $\mathcal{C}_k^-$ the sets of unbalanced signed complete graphs and of negative cycles of order $k$, respectively. 
	Below we list two results in their paper and provide a commentary. Note that, in Figure \ref{fig:ExtreGs}, black solid lines (resp., dashed lines) represent positive edges (resp., negative edges), golden lines between two vertex sets represent the connection of all possible positive edges, and the solid circle represent all-positive signed complete graph.
	\begin{theorem}\label{tho:edge-c_3free}\cite[Theorem 1.2]{wang2022extremed}
		If $\Gamma$ is a connected $\mathcal{C}_3^-$-free unbalanced signed graph of order $n$, then 
		$$e(\Gamma) \leq \frac{n(n - 1)}{2} - (n - 2), $$
		with equality holding if and only if $\Gamma \sim \Gamma^{s,t}$ (see Fig. \ref{fig:ExtreGs}), where $s + t = n -2$ and $s,t \geq 1.$
	\end{theorem}
	\begin{theorem}\label{tho:spectra-c_3free}\cite[Theorem 1.3]{wang2022extremed}
		If $\Gamma$ is a connected $\mathcal{C}_3^-$-free unbalanced signed graph of order $n$, then 
		$$\rho(\Gamma) \leq \frac{1}{2}(\sqrt{n^2 - 8} + n - 4), $$
		with equality holding if and only if $\Gamma \sim \Gamma^{1,n-3}.$
	\end{theorem}	
	
	\begin{remark}
		In Theorems \ref{tho:edge-c_3free} and \ref{tho:spectra-c_3free}, the conditions of connectivity can be omitted. 
		In fact, for the former theorem, if the signed edge extremal graph is not connected, we can add an edge between two connected components to obtain a $\mathcal{C}^-_{3}$-free unbalanced signed graph, leading to a contradiction. 
		For the latter theorem, we know that for a $\mathcal{C}^-_{3}$-free unbalanced signed graph with the maximum spectral radius, its spectral radius is equal to its index. 
		Therefore, according to Proposition \ref{pro:perturbation}, which we will prove later, if the signed spectral extremal graph is not connected, we can add a positive edge between two connected components, also leading to a contradiction.
	\end{remark}
	Recently, several researches have explored related issues (see \cite{chen2023turan,wang2023largest,wang2023spectral,wang2023turan}), and here we recall several results that will be helpful for study later.
	Chen and Yuan, in \cite{chen2023turan}, investigated the signed edge Tur{\'a}n problem and signed spectral Tur{\'a}n problem for $\mathcal{K}_4^-$-free signed graphs, and their results are as follows: 
	\begin{theorem}\label{tho:edge-k_4-free}\cite[Theorem 1.5]{chen2023turan}
		If $\Gamma$ is a $\mathcal{K}_4^-$-free unbalanced signed graph of order $n$  $ (n \geq 7)$.
		Then $$e(\Gamma) \leq \frac{n(n - 1)}{2} - (n - 3).$$
	\end{theorem}
	
	Note that by consulting the tables of signed graphs with order at most 6\cite{bussemaker1991tables}, the condition $n \geq 7$ can be removed from the theorem above. 
	\begin{theorem}\label{tho:spectra-k_4-free}\cite[Theorem 1.6]{chen2023turan}
		If $\Gamma$ is a $\mathcal{K}_4^-$-free unbalanced signed graph of order $n$.
		Then $$\rho(\Gamma) \leq n - 2,$$
		with equality holding if and only if $\Gamma \sim \Gamma_{1,1}.$ (see Fig. \ref{fig:ExtreGs})
	\end{theorem}
  Wang \cite{wang2023spectral} continued such studies on $\mathcal{K}_5^-$-free signed graph and left an open problem: ``What is the maximum spectral radius among all $\mathcal{K}_{r+1}^-$-free ($r \geq 5$) unbalanced signed graphs of order $n$?"
  In this work, we answer above problem for $3 \leq r \leq \lfloor \frac{n}{2} \rfloor$. 
  Before getting this result, we give a tight upper bound on the number of edges of $\mathcal{K}_{r + 1}^-$-free unbalanced signed graphs.
  
  \begin{theorem}\label{tho:edge-k_{r+1}-free}
  	If $\Gamma$ is a $\mathcal{K}_{r+1}^-$-free unbalanced signed graph of order $n$.
  	Then $$e(\Gamma) \leq \frac{n(n - 1)}{2} - (n - r).$$
  \end{theorem}
  
  Observe that $e(\Gamma_{1,r-2}) = \frac{n(n - 1)}{2} - (n - r)$, so the upper bound above is a supremum. 
  Below we show that the $\mathcal{K}^-_{r + 1}$-free unbalanced signed graph with maximum index is switching equivalent to $\Gamma_{1,r-2}$.
  
  \begin{theorem}\label{tho:spectra-k_{r+1}-free}
  	If $\Gamma$ is a $\mathcal{K}_{r+1}^-$-free $(r \geq 3)$ unbalanced signed graph of order $n$.
  	Then $$\lambda_{1}(\Gamma) \leq \lambda_1(\Gamma_{1,r-2}),$$
  	with equality holding if and only if $\Gamma \sim \Gamma_{1,r-2}$ (see Fig. \ref{fig:ExtreGs}).
  \end{theorem}
	
	The proposition below add some numerical estimates to Theorem \ref{tho:spectra-k_{r+1}-free}.
	\begin{proposition}\label{pro:num}
		The index $\lambda_{1}(\Gamma_{1,r-2})$ corresponds to the largest root of the polynomial $$f(x) = x^3 + (3-n)x^2 + (3-n-r)x + (n + 4)r - (r^2 + n + 7),$$ and satisfies
		$$n - 2 \leq \lambda_{1}(\Gamma_{1,r-2}) < n - 1.$$
	\end{proposition}
	
	The remainder of this article is organized as follows: in Section \ref{sec:pre}, we present some fundamental properties and conclusions that will be used in the sequel. 
	In Section \ref{sec:pro}, we provide the proofs of Theorems \ref{tho:edge-k_{r+1}-free}, \ref{tho:spectra-k_{r+1}-free}, and Proposition \ref{pro:num}. 
	In the concluding remarks, we analyze our results and provide some comments.

\indent

\section{Preliminaries}\label{sec:pre}
  In this section we introduce some results which will be useful in the sequel.
  The lemma below concerns equitable partition. 
  Consider a partition $P = \{V_1, \cdots, V_m\}$ of the set $V = \{1, \cdots, n\}.$ 
  The characteristic matrix $\chi_P$ of $P$ is the $n \times m$ matrix whose columns are the characteristic vectors of $ V_1, \cdots, V_m.$
  Consider a symmetric matrix $M$ of order $n$, with rows and columns are partitioned according to $P$. 
  The partition of $M$ is equitable if each submatrix $M_{i,j}$ is formed by the rows of $V_i$ and the columns of $V_j$ has
  constant row sums $q_{i,j}$. 
  The $m \times m$ matrix $Q = (q_{i,j} )_{1 \leq i,j \leq m} $ is called the quotient matrix of $M$ with respect to the equitable partition $P$.
  \begin{lemma}\cite[p. 30]{brouwer2022spectra}\label{lem:quotient}
  	The matrix $M$ has the following two kinds of eigenvectors and eigenvalues:\\
    \indent (\romannumeral1) The eigenvectors in the column space of $\chi_P$; the corresponding eigenvalues coincide with the eigenvalues of $Q$,\\
  	\indent (\romannumeral2) The eigenvectors orthogonal to the columns of $\chi_P$; the corresponding eigenvalues
  	of $M$ remain unchanged if some scalar multiple of the all-one block $J$ is added to block
  	$M_{i,j}$ for each $i, j \in \{1, \cdots, m\}.$
  \end{lemma}

   Next we present a celebrated result of signed graphs, which is an important method for determining the switching equivalence of two signed graphs with the same underlying graph. 
  \begin{lemma}\cite[Proposition 3.2]{zaslavskyi1982signed}\label{lem:swi}
  	Two signed graphs on the same underlying graph are switching  equivalent if and only if they have the same list of balanced cycles. 
  \end{lemma}
  
  The clique number of a graph $G$, denoted by $\omega(G)$, is the maximum order of a clique in $G$. The balanced clique number of a signed graph $\Gamma$, denoted by $\omega_b(\Gamma)$, is the maximum order of a balanced clique in $\Gamma$.
  The following two lemmas give two upper bounds for the index of a graph (resp., a signed graph) in terms of the order and clique number (resp., balanced clique number).
  
  \begin{lemma}\cite{wilf1986spectral}
  	Let $G$ be a graph of order $n$. Then
  \begin{equation*}
  	\lambda_{1}(G) \leq n(1-\frac{1}{\omega(G)}).
  \end{equation*}
  \end{lemma}
  
   \begin{lemma}\cite[Prosition 5]{wang2021eigenvalues}\label{lem:bal}
  	Let $\Gamma$ be a signed graph of order $n$. Then
  	\begin{equation*}
  		\lambda_{1}(\Gamma) \leq n(1-\frac{1}{\omega_b(\Gamma)}).
  	\end{equation*}
  \end{lemma}

  We end this section by following lemma which says that the index of a signed graph is not larger than that of its underlying graph. 
  
  \begin{lemma}\cite[Theorem 2.1 and Proposition 2.5]{brunetti2022unbalanced}\label{lem:rel}
  	For a non-empty signed graph $\Gamma$ of order $n$, $\lambda_{1}(\Gamma) \leq \lambda_{1}(G)$.
	Furthermore, $\lambda_1(\Gamma) \leq n - 1$, with equality if and only if $\Gamma$ is balanced and complete.
  \end{lemma}

\section{ Proofs}\label{sec:pro}
	This section is devoted to the proofs of Theorem \ref{tho:edge-k_{r+1}-free}, \ref{tho:spectra-k_{r+1}-free}, and Proposition \ref{pro:num}.
	For signed graph notation and concepts undefined here, we refer the reader to \cite{zaslavskyi1982signed}. For introductory material on the signed graph theory see the survey of Zaslavsky \cite{zaslavskyi2008matrices}
	and its references. 
	In particular, let $\Gamma$ be a signed graph, and $X$ and $Y$ be disjoint sets of vertices of $\Gamma$.
	We write:\\
	\indent - $V(\Gamma) = \{v_1, \cdots, v_n\}$ for the set of vertices of $\Gamma$, and $e(\Gamma)$ for the number of its edge;\\
	\indent	- $\Gamma[X]$ for the signed graph induced by $X$, and $e(X)$ for $e(\Gamma[X])$;\\
	\indent	- $e(X,Y)$ for the number of edges joining vertices in $X$ to vertices in $Y$;\\
	\indent - $N_\Gamma(v)$ for the set of neighbors of a vertex $v$ in $\Gamma$, and $d_\Gamma(v)$ for $|N_\Gamma(v)|$.
	
		\begin{proof}[\indent Proof of Theorem \ref{tho:edge-k_{r+1}-free}]
		From Theorem \ref{tho:edge-k_4-free} we know that the theorem holds for $r = 3$. 
		We will prove our result by induction on the order of the forbidden unbalanced signed complete graphs.
		Assume that the theorem is true for all values not greater than $r$ and we prove it for $r + 1$.
		Let $\Gamma$ be a $\mathcal{K}^-_{r+1}$-free unbalanced signed graph with maximum possible number of edges.
		Note that $\Gamma_{1,r-2}$ is $\mathcal{K}_{r+1}^-$-free, and so $e(\Gamma) \geq e(\Gamma_{1,r-2}).$\\
		\indent First we claim that $\Gamma$ contains an unbalanced $K_r$. 
		Otherwise, by our hypothesis we know 
		$$e(\Gamma) \leq \frac{n(n-1)}{2} -  (n - r + 1) < e(\Gamma_{1,r-2}),$$ a contradiction. 
		Let $X$ be the vertex set of an unbalanced signed complete graph with order $r$ and let $Y$ be its complement. Since each vertex in $Y$ can have at most $r - 1$ neighbours in $X$, the number of edges between $X$ and $Y$ is at most $(r - 1)(n - r)$.
		We see that
		$$ e (\Gamma) = e(X) + e(Y) + e(X,Y) \leq \tbinom{r}{2} + \tbinom{n - r}{2} + (r - 1)(n - r) = \frac{n(n - 1)}{2} - (n - r).$$ The theorem follows.
 	\end{proof}
	Note that the approach for proving Theorem \ref{tho:edge-k_{r+1}-free} is inspired by the proof of Tur{\'a}n's Theorem, with the key difference being that we use induction on $r$ here, whereas his proof employs induction on $n$.
 
	\begin{proof}[\indent Proof of Proposition \ref{pro:num}]
		We give a vertex partition as $V_1 = \{v_1 \}$, $V_2 = \{v_2\}$, $V_{3} = \{v_3, \cdots, v_r\}$, and $V_4 = \{v_{r+1}, \cdots, v_n\}$. Then the adjacency matrix of $\Gamma_{1,r-2}$ and its quotient matrix $Q$ are
		\[A(\Gamma_{1,r-2}) = \bbordermatrix{
			& V_1  & V_2 & V_3 & V_4 & \cr
			V_1& 0 & -1 & \textbf{j}^{\top} & \textbf{0}^{\top} \cr
			V_2 & -1 & 0 & \textbf{j}^{\top} & \textbf{j}^{\top} \cr
			V_3 & \textbf{j} & \textbf{j} & J-I & J \cr
			V_4& \textbf{0} & \textbf{j} & J & J-I \cr
		} and ~
		Q_1 = \bbordermatrix{
			& V_1  & V_2 & V_3 & V_4 & \cr
			V_1& 0 & -1 & r-2 & 0 \cr
			V_2 & -1 & 0 & r-2 & n-r \cr
			V_3 & 1 & 1 & r-3 & n-r \cr
			V_4& 0 & 1 & r-2 & n-r-1  \cr
		},
		\] where \textbf{0} and \textbf{j} represent the zero vector and the all-ones vector of appropriate dimensions, respectively, and $I$ and $J$ denote the identity matrix and all-ones matrix of appropriate orders, respectively. By Lemma \ref{lem:quotient}, the eigenvalues of $Q_1$ are that of $A(\Gamma_{1,r-2})$ and the other eigenvalues of $A(\Gamma_{1,r-2})$ remain if we add some scalar multiple of $J$ from the blocks equal to $-1$, $\textbf{j}$, $J$, and $J - I$.
		Then $A(\Gamma_{1,r-2})$ and $Q_1$ become
		\[A^\prime(\Gamma_{1,r-2}) = \bbordermatrix{
			& V_1  & V_2 & V_3 & V_4 & \cr
			V_1& 0 & 0 & \textbf{0}^{\top} & \textbf{0}^{\top} \cr
			V_2 & 0 & 0 & \textbf{0}^{\top} & \textbf{0}^{\top} \cr
			V_3 & \textbf{0} & \textbf{0} & -I & 0 \cr
			V_4& \textbf{0} & \textbf{0} & 0 & -I \cr
		} and ~
		Q_1^\prime = \bbordermatrix{
			& V_1  & V_2 & V_3 & V_4 & \cr
			V_1& 0 & 0 & 0 & 0 \cr
			V_2 & 0 & 0 & 0 & 0 \cr
			V_3 & 0 & 0 & -1 & 0 \cr
			V_4& 0 & 0 & 0 & -1  \cr
		},
		\]
		The eigenvalues of matrix $A^\prime(\Gamma_{1,r-2})$ except the eigenvalues of $Q^\prime$ are $-1$ with multiplicity $n-4$.
		Then the eigenvalues of $\Gamma_{1,r-2}$ are the eigenvalues of $Q_1$ and $-1$ with multiplicity $n-4$.
		Therefore, $\lambda_{1}(\Gamma_{1,r-2}) = \lambda_1(Q_1)$.
		By direct calculation, the characteristic polynomial of the matrix $Q_1$ is 
		$$g(x) = (x+1) (x^3 + (3-n)x^2 + (3-n-r)x + (n + 4)r - (r^2 + n + 7)) = (x+1)f(x).$$
		Thus $\lambda_{1}(\Gamma_{1,r-2})$ is the largest root of $f(x) = x^3 + (3-n)x^2 + (3-n-r)x + (n + 4)r - (r^2 + n + 7) = 0$.
		By simple calculations $f(n-2) = - (r - 3) \leq 0$. 
		So we have $\lambda_{1}(\Gamma_{1,r-2}) \geq n-2$, and $\lambda_{1}(\Gamma_{1,r-2})<n-1$ form Lemma \ref{lem:rel}.
	\end{proof}
	
	To simplify the proof of Theorem \ref{tho:spectra-k_{r+1}-free}, we shall prove several auxiliary statements.
	First, note that according to Perron-Frobenius theory, there exists a strictly positive eigenvector corresponding to the index of a simple connected graph $G$. 
	Furthermore, if we add some edges in $G$, the index of the resulting graph is larger than that of $G$. 
	These are usually incorrect for signed graphs, but we have following results.
	
	\begin{lemma}\cite[Lemma1]{stanic2018perturbations}\label{lem:nonneg}
		Let $\Gamma$ be a signed graph with $n$ vertices. 
		Then there exists a signed graph $\Gamma^\prime$ switching equivalent to $\Gamma$ such that $\lambda_1(\Gamma^\prime)$ has a non-negative eigenvector.
	\end{lemma}

Let $x = (x_1, x_2, \cdots, x_n)^{\intercal}$ be an eigenvector corresponding to the index $\lambda_{1}(\Gamma)$ of a signed graph $\Gamma$. 
The entry $x_i$ is usually corresponding to the vertex $v_i$ of $\Gamma$. 
So the eigenvalue equation for $v_i$ reads as follows
\begin{align*}
	\lambda_{1}(\Gamma)x_{i}  = & \sum_{v_j \in N_{\Gamma}(v_i)} \sigma(v_iv_j)x_j. 
\end{align*}
The following lemma can be proved based on \cite[Theorem 3]{stanic2018perturbations} or \cite[Proposition 2.1]{akbari2019largest}.
For the sake of completeness, here we provide a proof.

\begin{proposition}\label{pro:perturbation}
	Let $\Gamma = (G, \sigma)$ be a signed graph with a non-negative unit eigenvector $x = (x_1, x_2, \cdots, x_n)^{\intercal}$ corresponding to the largest eigenvalue $\lambda_{1}(\Gamma)$. 
	If we perform one of the following perturbations in $\Gamma$:\\
	\indent \textbf{(\rmnum{1})} Adding some positive edges,\\
	\indent \textbf{(\rmnum{2})} Removing some negative edges,\\
	\indent \textbf{(\rmnum{3})} Reversing the signs of some negative edges,\\
	resulting in a new signed graph $\Gamma ^ \prime $, then $\lambda_{1}(\Gamma^\prime) \geq \lambda_{1}(\Gamma)$. 
	The equality holds if and only if the entries of $x$  corresponding to the endpoints of these edges are all zeros.\\
	\indent And if we perform one of the following perturbations in $\Gamma$:\\
	\indent \textbf{(\rmnum{4})} Rotating the positive edge $v_iv_j$ to the non-edge position $v_iv_k$, where $x_j \leq x_k$,\\
	\indent \textbf{(\rmnum{5})} Reversing the sign
	of the positive edge $v_iv_j$ and the negative edge $v_iv_k$, where $x_j \leq x_k$,\\
	resulting in a new signed graph $\Gamma ^ \prime $, then $\lambda_{1}(\Gamma^\prime) \geq \lambda_{1}(\Gamma)$. 
	The equality holds if and only if $x_i = 0$ and $x_j = x_k$.\\
\end{proposition}	
\begin{proof}
	For \textbf{(\rmnum{1})}, we denote by $E_1 \subseteq E(\Gamma^\prime)$ the set of added positive edges.
	By Rayleigh Principle, we have
	\begin{align*}
		\lambda_1({\Gamma^\prime}) - \lambda_1(\Gamma) 
		& =  \max \limits_{\Vert y \Vert = 1} y^\intercal A(\Gamma^\prime) y - x^\intercal A(\Gamma) x\\
		& \geq x^\intercal A(\Gamma^\prime) x - x^\intercal A(\Gamma) x\\
		& = 2 \sum_{v_iv_j \in E_1 } x_{i}x_{j}\\
		& \geq 0.
	\end{align*}
	If $\lambda_{1}(\Gamma^\prime) = \lambda_{1}(\Gamma)$, then all the equalities hold and so $x$ is an eigenvector of $A(\Gamma^\prime)$ corresponding to the eigenvalue $\lambda_1({\Gamma^\prime})$. 
	Take one positive edge from $E_1$, say $v_kv_l$.
	We will show $x_l = 0$.  
	Assume that the added positive edges with one endpoint $v_k$ are $v_kv_l, v_kv_{k_1}, v_kv_{k_2}, \cdots, v_kv_{k_s}$.
	According to the following eigenvalue equations,
	\begin{align*}
		\lambda_{1}(\Gamma)x_{k} & = \sum_{v_h \in N_{\Gamma}(v_{k})} \sigma(v_hv_{k})x_h,\\
		\lambda_{1}(\Gamma^\prime)x_{k} & = \sum_{v_h \in N_{\Gamma}(v_{k})} \sigma(v_hv_{k})x_h + \sum_{j = 1}^{s} x_{k_j} + x_l,\\
	\end{align*}
	we obtain $x_l = x_{k_1} = \cdots = x_{k_j} = 0$. 
	By similar analysis as above, the entries of $x$ corresponding to the endpoints of added positive edges are all zeros. \\
	\indent For \textbf{(\rmnum{2})}, we denote by $E_2 \subseteq E(\Gamma)$ the set of deleted negative edges.
	We have
	\begin{align*}
		\lambda_1({\Gamma^\prime}) - \lambda_1(\Gamma) 
		& \geq x^\intercal A(\Gamma^\prime) x - x^\intercal A(\Gamma) x\\
		& = 2 \sum_{v_iv_j \in E_2 } x_{i}x_{j}\\
		& \geq 0.
	\end{align*}
	If $\lambda_{1}(\Gamma^\prime) = \lambda_{1}(\Gamma)$, then all the equalities hold and so $x$ is an eigenvector of $A(\Gamma^\prime)$ corresponding to the eigenvalue $\lambda_1({\Gamma^\prime})$. 
	Take one negative edge from $E_2$, say $v_kv_l$.
	We will show $x_l = 0$.  
	Assume that the deleted negative edges with one endpoint $v_k$ are $v_kv_l, v_kv_{k_1}, v_kv_{k_2}, \cdots, v_kv_{k_s}$.
	According to the following eigenvalue equations,
	\begin{align*}
		\lambda_{1}(\Gamma)x_{k} & = \sum_{v_h \in N_{\Gamma}(v_{k})} \sigma(v_hv_{k})x_h -  \sum_{j = 1}^{s} x_{k_j} - x_l,\\
		\lambda_{1}(\Gamma^\prime)x_{k} & = \sum_{v_h \in N_{\Gamma}(v_{k})} \sigma(v_hv_{k})x_h,\\
	\end{align*}
	we obtain $x_l = x_{k_1} = \cdots = x_{k_j} = 0$. 
	By similar analysis as above, the entries of $x$ corresponding to the endpoints of deleted negative edges are all zeros.\\
	\indent For \textbf{(\rmnum{3})}, we denote by $E_3 \subseteq E(\Gamma)$ the set of changed negative edges.
	We have
	\begin{align*}
		\lambda_1({\Gamma^\prime}) - \lambda_1(\Gamma) 
		& \geq x^\intercal A(\Gamma^\prime) x - x^\intercal A(\Gamma) x\\
		& = 4 \sum_{v_iv_j \in E_3 } x_{i}x_{j}\\
		& \geq 0.
	\end{align*}
	If $\lambda_{1}(\Gamma^\prime) = \lambda_{1}(\Gamma)$, then all the equalities hold and so $x$ is an eigenvector of $A(\Gamma^\prime)$ corresponding to the eigenvalue $\lambda_1({\Gamma^\prime})$. 
	Take one negative edge from $E_3$, say $v_kv_l$.
	We will show $x_l = 0$.  
	Assume that the changed negative edges with one endpoint $v_k$ are $v_kv_l, v_kv_{k_1}, v_kv_{k_2}, \cdots, v_kv_{k_s}$.
	According to the following eigenvalue equations,
	\begin{align*}
		\lambda_{1}(\Gamma)x_{k} & = \sum_{v_h \in N_{\Gamma}(v_{k})} \sigma(v_hv_{k})x_h -  \sum_{j = 1}^{s} x_{k_j} - x_l,\\
		\lambda_{1}(\Gamma^\prime)x_{k} & = \sum_{v_h \in N_{\Gamma}(v_{k})} \sigma(v_hv_{k})x_h + \sum_{j = 1}^{s} x_{k_j} + x_l,\\
	\end{align*}
	we obtain $x_l = x_{k_1} = \cdots = x_{k_j} = 0$. 
	By similar analysis as above, the entries of $x$ corresponding to the endpoints of changed negative edges are all zeros.\\
	\indent For \textbf{(\rmnum{4})}, we have
	\begin{align*}
		\lambda_1({\Gamma^\prime}) - \lambda_1(\Gamma) 
		& \geq x^\intercal A(\Gamma^\prime) x - x^\intercal A(\Gamma) x\\
		& = 2x_{i}( x_{k} - x_{j})\\
		& \geq 0.
	\end{align*} 
	If $\lambda_{1}(\Gamma^\prime) = \lambda_{1}(\Gamma)$, then all the equalities hold and so $x$ is an eigenvector of $A(\Gamma^\prime)$ corresponding to the eigenvalue $\lambda_1({\Gamma^\prime})$. 
	In view of the following eigenvalue equations,
		\begin{align*}
		\lambda_{1}(\Gamma)x_{i} & = \sum_{v_h \in N_{\Gamma}(v_{i}) \setminus v_j}  \sigma(v_hv_{i})x_h + x_j,\\
		\lambda_{1}(\Gamma^\prime)x_{i} & = \sum_{v_h \in N_{\Gamma}(v_{i}) \setminus v_j}  \sigma(v_hv_{i})x_h + x_k,\\
		\lambda_{1}(\Gamma)x_{j} & = \sum_{v_h \in N_{\Gamma}(v_{j}) \setminus v_i}  \sigma(v_hv_{j})x_h + x_i,\\
		\lambda_{1}(\Gamma^\prime)x_{j} & = \sum_{v_h \in N_{\Gamma}(v_{j}) \setminus v_i}  \sigma(v_hv_{j})x_h,\\
		\lambda_{1}(\Gamma)x_{k} & = \sum_{v_h \in N_{\Gamma}(v_{k})}  \sigma(v_hv_{k})x_h,\\
		\lambda_{1}(\Gamma^\prime)x_{k} & = \sum_{v_h \in N_{\Gamma}(v_{k})}  \sigma(v_hv_{k})x_h + x_i,\\
	\end{align*}
	we have $x_i = 0$ and $x_j = x_k$.\\
	\indent For \textbf{(\rmnum{5})}, we have
	\begin{align*}
		\lambda_1({\Gamma^\prime}) - \lambda_1(\Gamma) 
		& \geq x^\intercal A(\Gamma^\prime) x - x^\intercal A(\Gamma) x\\
		& = 4x_{i}( x_{k} - x_{j})\\
		& \geq 0.
	\end{align*} 
	The remainder proof are almost the same as that of \textbf{(\rmnum{4})} and we omit it. The converse is clear.
\end{proof}

	\begin{lemma}\label{lem:zer}
	Let $\Gamma = (G, \sigma)$ be a signed graph with a unit eigenvector $x = (x_1, x_2, \cdots, x_n)^{\intercal}$ corresponding to $\lambda_{1}(\Gamma)$. 
	If $\lambda_1(\Gamma) > n - k$, then $x$ has at most $k - 2$ zero component.
	\begin{proof}
		Without loss of generality, assume for a contradiction that $x_1 = x_2 = \cdots = x_{k - 1} = 0 $.
		Deleting the corresponding vertices from $\Gamma$ to obtain a signed graph $\Gamma^\prime$.
		Then by Rayleigh Principle and Lemma \ref{lem:rel},
		\begin{align*}
			\lambda_{1}(\Gamma) &= (x_{k}, \cdots, x_n) A(\Gamma^\prime) (x_{k}, \cdots, x_n)^\intercal\\
			&\leq\lambda_{1} (\Gamma^\prime) \leq \lambda_{1}(K_{n - k + 1}) = n - k,						 
		\end{align*}
		a contradiction.
	\end{proof}
\end{lemma}

\begin{remark}\label{rem:connectedandpositive}
	Note from Proposition \ref{pro:num} that $\Gamma_{1,r-2}$ is a $\mathcal{K}_{r+1}^-$-free unbalanced signed graph with index $\lambda_{1}(\Gamma_{1,r-2}) > n - 2$, for $r \geq 4$. 
	Combining this with Proposition \ref{pro:perturbation} and Lemma \ref{lem:zer},
	we know that if $\Gamma$ is a $\mathcal{K}_{r+1}^-$-free unbalanced signed graph with maximum index, then it is connected.
	Furthermore, if $r \geq 4$ we can find an eigenvector corresponding to $\lambda_{1}(\Gamma)$ with no zero components.
\end{remark}

\begin{proof}[\indent Proof of Theorem \ref{tho:spectra-k_{r+1}-free}]
	From Theorem \ref{tho:spectra-k_4-free} we know that our theorem holds when $r = 3$.
	By induction on $r$, assume it is true for all the values not greater than $r$ and prove it for $r+1$.
	Assume that $\Gamma$ has the maximum index over all $\mathcal{K}_{r+1}^-$-free unbalanced signed graphs. 
	In view of Lemma \ref{lem:nonneg} and Remark \ref{rem:connectedandpositive} we can find $\widetilde{\Gamma} \sim \Gamma$ with a positive eigenvector $x = (x_1, \cdots, x_n)^{\top}$ corresponding to $\lambda_1(\widetilde{\Gamma}) > n - 2$.
	Note from Lemma \ref{lem:swi} that $\widetilde{\Gamma}$ is also a  connected $\mathcal{K}^-_{r+1}$-free unbalanced signed graph. We will show $\widetilde{\Gamma} = \Gamma_{1,r-2}$ step by step.\\
	\indent First, since $\widetilde{\Gamma}$ is unbalanced, there exist at least one negative edge and at least one negative cycle. 
	Take a negative cycle $\mathcal{C} = v_1v_2 \cdots v_lv_1$ of the shortest length from $\widetilde{\Gamma}$. 
	We claim $l = 3$, otherwise $\widetilde{\Gamma}$ is $\mathcal{C}_3^-$-free, and so $\lambda_{1}(\widetilde{\Gamma}) \leq \frac{1}{2}(\sqrt{n^2 - 8} + n - 4) < n - 2$ by Theorem \ref{tho:spectra-c_3free}, a contraction.
	Thus, $\mathcal{C}$ is a negative triangle on vertices $v_1$, $v_2$, and $v_3$.\\
	\indent Secondly, we say that all the negative edges of $\widetilde{\Gamma}$ are contained in $\mathcal{C}$. Indeed, if there exists a negative edge not in $\mathcal{C}$, by Proposition \ref{pro:perturbation} we may delete it resulting a $\mathcal{K}_{r+1}^-$-free unbalanced signed graph with larger	 index, a contradiction.
	Therefore, the number of negative edges of $\widetilde{\Gamma}$ is either $1$ or $3$.\\
	\indent We conclude that $\widetilde{\Gamma}$ contains only one negative edge. 
	Actually, if not, then $\mathcal{C}$ is a negative triangle with three negative edges and by Proposition \ref{pro:perturbation} we may reverse signs of two of those, resulting a $\mathcal{K}_{r+1}^-$-free unbalanced signed graph with larger index, which leads a contradiction.\\
	\indent Next we claim that $\widetilde{\Gamma}$ contains an unbalanced $K_r$, written as $K_r^-$. 
	If not, $\widetilde{\Gamma}$ is $\mathcal{K}_r^-$-free and by induction hypothesis we know that $\widetilde{\Gamma} \sim \Gamma_{1,r-3}$ which has smaller index than $\Gamma_{1,r-2}$ by Proposition $\ref{pro:perturbation}$, a contradiction. 
	Without loss of generality, suppose that $X = V(K_r^-) = \{v_1, \cdots, v_r\}$ and $Y = V(\widetilde{\Gamma}) \setminus X$, and further that $x_1 \leq x_2$ and $x_3 \leq \cdots \leq x_r$.
	Let $W_1 = N_{\widetilde{\Gamma}}(v_1) \setminus N_{\widetilde{\Gamma}}(v_2),	W_2 = N_{\widetilde{\Gamma}}(v_2) \setminus N_{\widetilde{\Gamma}}(v_1)$, and $W = N_{\widetilde{\Gamma}}(v_1) \cap N_{\widetilde{\Gamma}}(v_2) \setminus X $.
	We claim that $W_1 = \emptyset$, otherwise there exist a vertex $v_k$ satisfies $v_k \stackrel{+}{\sim} v_1$ and $v_k \not \sim v_2$,
	and by Proposition \ref{pro:perturbation}  we can rotate the positive edge $v_kv_1$ to the non-edge position $v_kv_2$, getting a $\mathcal{K}_{r+1}^-$-free unbalanced  signed graph with larger index than $\widetilde{\Gamma}$, a contradiction. \\
	\indent We proceed with our proof and establish that $x_2 < x_3$. 
	Assume for a contradiction that $x_2 \geq x_3$. 
	If there exists $V_1 \subseteq V(\widetilde{\Gamma})$ such that $\widetilde{\Gamma}[V_1 \cup \{v_1, v_3\}]$ is $(K_{r+1},+)$, then by $W_1 = \emptyset$ we have that each vertex in $V_1$ is adjacent to $v_2$ and so $\widetilde{\Gamma}[V_1 \cup \{v_1, v_2\}]$ is an unbalanced $K_{r+1}$, a contradiction.
	So there exist no $V_1 \subseteq V(\widetilde{\Gamma})$ such that $\widetilde{\Gamma}[V_1 \cup \{v_1, v_3\}]$ is $(K_{r+1},+)$, and  by Proposition \ref{pro:perturbation}  we may reverse the signs of $v_1v_2$ and $v_1v_3$, resulting a $\mathcal{K}_{r+1}^-$-free  unbalanced signed graph with larger index than $\widetilde{\Gamma}$, a contradiction.
	Then we know that $x_1 \leq x_2 < x_3 \leq \cdots \leq x_r$.\\
	\indent Note that each vertex in $Y$ is adjacent to at most $r-1$ vertices in $X$.
	Then we claim that $W = \emptyset$. 
	Otherwise, there exists a vertex $v_i \in W$ not being adjacent to a vertex $v_j \in X \setminus \{v_1, v_2\}$, and then we can rotate the positive edge $v_iv_1$ to the non-edge position $v_iv_j$.
	By Proposition \ref{pro:perturbation} and the fact $x_j > x_1$, a contradiction.\\
	\indent Summing up, $\widetilde{\Gamma}$ must be a subgraph of $\Gamma_{1,r-2}$ and according to Proposition \ref{pro:perturbation} it actually is $\Gamma_{1,r-2}$, completing the proof of Theorem \ref{tho:spectra-k_{r+1}-free}.
\end{proof}

\section{Concluding remarks}\label{sec:con}
	In Theorem \ref{tho:spectra-k_{r+1}-free}, the signed spectral extremal graph has the number of edges reaching the upper bound stated in Theorem \ref{tho:edge-k_{r+1}-free}. 
	Hence, we can conclude that, among all $\mathcal{K}_{r+1}^-$-free $(r \geq 3)$ unbalanced signed graph, the signed spectral extremal graph must be a signed edge extremal graph. 
	Unfortunately, the signed edge extremal graph is not unique, up to switching equivalence. 
	In fact, there are various possible signed graphs that can attain the upper bound on the number of edges, as discussed in \cite[Theorem 1.5]{chen2023turan} for the case of $r = 3$.
	
	In view of the following remark, we show that Theorem \ref{tho:spectra-k_{r+1}-free} partly solve the problem ``What is the maximum spectral radius among all $\mathcal{K}_{r+1}^-$-free ($r \geq 5$) unbalanced signed graphs of order $n$?"
	The negation of $\Gamma$ (denoted by  $-\Gamma$) is obtained by reversing the sign of each edge in $\Gamma$. 
	Clearly, the eigenvalues of $-\Gamma$ are obtained by reversing the signs of the eigenvalues of $\Gamma$.
	\begin{remark}\label{rem:speradius}
		Let $\Gamma = (G, \sigma)$ be a signed graph of order $n \geq 2r$ for an integer $r \geq 4$.
		If $\rho(\Gamma) > n - 2$ and $ - \Gamma$ is $\mathcal{K}_{r+1}^+$-free, then $\rho(\Gamma) = \lambda_{1}(\Gamma)$.
	\end{remark}
	\begin{proof}
		Assume for a contradiction that $\rho(\Gamma) = -\lambda_{n}(\Gamma)$.
		Since $-\Gamma$ is $\mathcal{K}_{r+1}^+$-free,
		using Lemma \ref{lem:bal}, we have
		\begin{align*}
			n - 2 < \rho(\Gamma) = -\lambda_{n}(\Gamma) &= \lambda_{1}(-\Gamma) \leq n \cdot (1 - \frac{1}{w_b(-\Gamma)}) \leq \frac{r - 1}{r}n,\\
		\end{align*}
		which is contradict to $n \geq 2r$.
	\end{proof}
	Indeed, according to the proof of Theorem \ref{tho:spectra-k_{r+1}-free} we know that the conditions of Remark \ref{rem:speradius} are satisfying by the $\mathcal{K}_{r+1}^-$-free $(r \geq 4)$ unbalanced signed graph having the maximum index, which means we solve above problem under the restriction $2r \leq n$. 
	The case of $2r \geq n$ is left and seems more challenging for further study.\\
	
\noindent \textbf{Declaration of competing interest}

The authors declare that there is no conflict of interest.

\vskip 0.6 true cm
\noindent {\textbf{Acknowledgments}}

This research is supported by the National Natural Science Foundation of China (No.11971164).


\baselineskip=0.25in


\begin{thebibliography}{99}

\bibitem{akbari2019largest} S. Akbari, F. Belardo,  F. Heydari, M. Maghasedi, M. Souri, On the largest eigenvalue of signed unicyclic graphs. Linear Algebra Appl. 581 (2019) 145--162.

\bibitem{bollobas1978extremal} B. Bollob{\'a}s, Extremal Graph Theory, Academic Press, London, 1978.

\bibitem{brouwer2022spectra} A. E. Brouwer, W. H. Haemers, Spectra of Graphs, Springer, 2011.

\bibitem{brunetti2022unbalanced} M. Brunetti, Z. Stani{\'c}, Unbalanced signed graphs with extremal spectral radius or index, Comput. Appl. Math. 41(3) (2022): 118.

\bibitem{bussemaker1991tables} F. C. Bussemaker, P. J. Cameron, J. J. Seidel, S. V. Tsaranov, Tables of signed graphs, Eut Report 91-WSK-01, Eindhoven, 1991.

\bibitem{chen2023turan} F. Chen, X. Y. Yuan, Tur\'{a}n problem for $\mathcal{K}_4^-$-free signed graphs, arXiv preprint arXiv:2306.06655 (2023).

\bibitem{li2022survey} Y. T. Li, W. J. Liu, L. H. Feng, A survey on spectral conditions for some extremal graph problems, Adv. Math. (China) 51 (2) (2022) 193--258.

\bibitem{nikiforovi2007bounds} V. Nikiforov, Bounds on graph eigenvalues II, Linear Algebra Appl. 427 (2007) 183--189.

\bibitem{nikiforovi2011some} V. Nikiforov, Some new results in extremal graph theory, in: Surveys in Combinatorics 2011, in: London Math. Society Lecture Note Ser., vol. 392, 2011, pp. 141--181.

\bibitem{stanic2018perturbations} Z. Stani{\'c}, Perturbations in a signed graph and its index, Discuss. Math. Graph T. 38 (3) (2018): 841--852. 

\bibitem{turan1941an} P. Tur{\'a}n, On an extremal problem in graph theory, Mat. Fiz. Lapok (in Hungarian) 48 (1941) 436--452.

\bibitem{turan1954the} P. Tur{\'a}n, On the theory of graphs, Colloq. Math. 3 (1954) 19--30.

\bibitem{wang2022extremed} D. J. Wang, Y. P. Hou, D. Q. Li,  Extremal results for $C_3^-$-free signed graphs, Linear Algebra Appl. 681 (2024) 47--65.

\bibitem{wang2023turan} J. J. Wang, Y. P. Hou, X. Y. Huang, Tur\'{a}n problem for $C_{2k+1}^-$-free signed graph, arXiv preprint arXiv:2310.11061 (2023).

\bibitem{wang2021eigenvalues} W. Wang, Z. D. Yan, J. G. Qian, Eigenvalues and chromatic number of a signed graph, Linear Algebra Appl. 619 (2021) 137--145.

\bibitem{wang2023spectral} Y. A. Wang, Spectral Tur\'{a}n problem for $\mathcal{K}_5^-$-free signed graphs, arXiv preprint arXiv:2309.15434 (2023).

\bibitem{wang2023largest} Y. A. Wang, H. Q. Lin, The largest eigenvalue of $\mathcal {C} _4^{-} $-free signed graphs, arXiv preprint arXiv:2309.04101 (2023).

\bibitem{wilf1986spectral} H. Wilf, Spectral bounds for the clique and independence numbers of graphs, J. Comb. Theory, Ser. B 40 (1986) 113--117.

\bibitem{zaslavskyi1982signed} T. Zaslavsky, Signed graphs, Discrete Appl. Math. 4 (1982) 47--74.

\bibitem{zaslavskyi2008matrices} T. Zaslavsky, Matrices in the theory of signed simple graphs, in: B.D. Acharya, G.O.H. Katona,
J. Nesetril (Eds.), Advances in Discrete Mathematics and Applications: Mysore, Ramanujan Mathematical Society, Mysore, 2010. 2008, pp. 207--229.





\end{thebibliography}
\end{document}